\documentclass[aps,pre,twocolumn,footinbib,showpacs]{revtex4-2}
\usepackage{graphicx}
\usepackage{bm}
\usepackage[urlcolor=blue]{hyperref}
\usepackage{color}
\usepackage[utf8]{inputenc}
\hypersetup{colorlinks,citecolor=blue,linkcolor=blue}
\usepackage{nicefrac}
\usepackage{ORCIDinREVTeX}
\usepackage{pdfpages} 
\makeatletter
\AtBeginDocument{\let\LS@rot\@undefined}
\makeatother

\definecolor{darkgreen}{rgb}{0.05, .65, 0.05}

\def\xmid{x_{\nicefrac{1}{2}}}

\begin{document}

\pacs{05.40.Fb, 02.50.Ey, 02.50.-r}

\title{Anatomy of an extreme event: What can we infer about the history of a heavy-tailed random walk?}

\author{Wesley W. Erickson}
\orcid{0000-0002-8457-8488}
\affiliation{Oregon Center for Optical, Molecular, and Quantum Science and Department of Physics, 1274 University of Oregon, Eugene, Oregon 97403-1274, USA}

\author{Daniel A. Steck}
\orcid{0000-0002-9120-7650}
\affiliation{Oregon Center for Optical, Molecular, and Quantum Science and Department of Physics, 1274 University of Oregon, Eugene, Oregon 97403-1274, USA}

\begin{abstract}
Extreme events are by nature rare and difficult to predict, yet are often
much more important than frequent, typical events.
An interesting counterpoint to the \textit{prediction} of such events is
their \textit{retrodiction}---given a process in an outlier state,
how did the events leading up to this endpoint unfold?
In particular, was there only a single, massive event, or was the history
a composite of multiple, smaller but still significant events?
To investigate this problem we take heavy-tailed stochastic 
processes (specifically,
the symmetric, $\alpha$-stable Lévy processes) as prototypical random walks.
A natural and useful
characteristic scale arises from the analysis of 
processes conditioned to arrive in a particular final state 
(Lévy bridges). For final displacements longer than this scale,
the scenario of a single, long jump is most likely, even though it
corresponds to a rare, extreme event.
On the other hand, for small final displacements, histories involving
extreme events tend to be suppressed.
To further illustrate the utility of this analysis,
we show how it provides an intuitive framework for understanding 
three problems related to boundary crossings of heavy-tailed processes. 
These examples illustrate how 
intuition fails to carry over from diffusive processes,
even very close to the Gaussian limit. One example
yields a computationally and conceptually useful representation of Lévy bridges
that illustrates how conditioning impacts the extreme-event content
of a random walk.
The other examples involve the conditioned boundary-crossing problem
and the ordinary first-escape problem; we discuss the observability
of the latter example in experiments with laser-cooled atoms.
\end{abstract}

\maketitle

\section{Introduction}

\begin{figure}[tb]
  \center
\includegraphics[width=\linewidth]{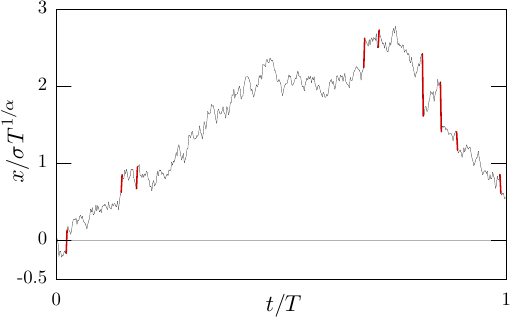}  \caption{
    Sample path of an unconditioned $\alpha$-stable 
    Lévy process, $\alpha=1.9$. Simulated path has
    time steps $\Delta t=T/500$ 
    with increments $\Delta x > L_\mathrm{b}(\Delta t/T)^{1/\alpha}$ 
    emphasized (bold/red).
  \label{fig:jumps_highlight}
}
\end{figure}

Extreme events affect us in many ways,
from geological and meteorological phenomena to 
market crashes and epidemics, and both science and society
have been increasingly appreciating the need to understand
and plan for such events \cite{albeverio06,beirlant2004}.
Gaussian stochastic models fail to predict extreme events, which
are commonly associated with probability distributions with
``heavy'' power-law tails.
Lévy processes (specifically, stable Lévy processes \cite{cont04, gardiner09, jacobs10})
in particular are important prototypes for
heavy-tailed random processes exhibiting large jumps or ``Lévy flights''
(Fig.~\ref{fig:jumps_highlight}), as they
are universal for random walks generated by
heavy-tailed distributions, in the same
sense that Gaussian processes are universal for finite-variance
steps.
Lévy processes
play an important role in understanding a wide range of phenomena
\cite{shlesinger95,uchaikin99}, including
ecology \cite{viswanathan96}, 
finance \cite{cont04},
fluid flows \cite{solomon93},
chaotic transport \cite{shlesinger93}, 
stochastic searches \cite{metzler09, palyulin14},
and particularly in laser-cooled atoms
\cite{marksteiner1996, katori1997, bardou02, sagi2012, kessler2012, barkai2014, afek2017, aghion2017}.
The stable processes also
produce strikingly counterintuitive behavior; for example,
intriguing work has shown that the image method fails to 
predict their first-passage times \cite{zumofen95,chechkin03,dybiec06,koren07}.

Of general importance in probability and statistics is the question of inference, 
which in stochastic processes is embodied by conditioned evolution.
The Brownian bridge---a continuous-time Gaussian stochastic
process specified to arrive at some final location (state)---is 
a well known and widely used examples of a conditioned process.
The properties and statistics of Brownian bridges
have been thoroughly studied \cite{borodin02};
they are productively applied in diverse areas, occurring in
financial mathematics \cite{brody07,moskowitz96},
models of animal movements \cite{horne07},
Monte Carlo methods in quantum mechanics \cite{gies03, mackrory16},
random interfaces and potentials \cite{levitz06, dean16, mori19},
and extreme-value statistics \cite{perret13}.
Because Lévy-type statistics arise in a similarly diverse
range of applications, and are also a cornerstone of extreme-event science,
clearly a detailed study of similarly conditioned, heavy-tailed processes is needed. 
 (An analogous generalization is to fractional Brownian bridges
 \cite{Delorme16}, which have been applied to the study
 of biological autoluminescence \cite{Dlask19}.)
Work on such Lévy bridges is at a nascent stage, however:
they have been formalized conceptually and applied to 
finance and insurance \cite{hoyle11,hoyle15}, and a few functionals of Lévy bridges
have been characterized \cite{fitzsimmons95,knight96,chaumont01}.

This paper explores the dynamics of 
continuous-time Lévy processes $x(t)$ conditioned to
arrive at the final state $L=x(T)$.
A key question that we address is: 
Was this arrival a result of a single, large event, or a composite of
multiple, smaller events?
From the typical behavior of heavy-tailed processes, where
rare but large events dominate the evolution,
one may expect that when arriving at an extreme state,
only a single extreme event is responsible, simply due to their rarity.
However, a proper accounting of the responsible events is
only possible by analyzing the conditional probabilities
for the state at intermediate times.
The structure of conditional probability densities 
for intermediate times $t\in(0,T)$ makes a transition from 
unimodal to bimodal as the arrival point $L$ varies,
leading to interesting and counterintuitive effects,
particularly in rare but important cases where an 
extreme jump occurred.
Above the bimodal transition, the typical conditioned history contains
only a single large event, while below
the transition the tendency is towards a composite of smaller events.
This analysis provides insight into first-passage 
problems for stable processes, highlighting
dramatic qualitative differences between Gaussian and
heavy-tailed processes, even when the latter are ``close to''
Gaussian. 
This work also provides a more precise, mathematical basis for the
intuition that random variations that occur in between rare, extreme events
tend to seem Gaussian, so much so that there is a strong temptation
to ignore extreme events in mathematical models, with sometimes devastating
consequences \cite{blackswan}.

A closely related existing result is the ``big-jump principle'' \cite{Vezzani19},
which observes under fairly general conditions
that for a sum of random variables, in the limit of a large
summed value, the distribution of the sum agrees with the distribution of the maximum
of the variables.  The implication is again that
extreme events are dominated by a single largest jump, rather than many small displacements.
This holds true even in the case of stretched-exponential processes, with
sub-power-law tails \cite{Burioni20}.
Another closely related concept is that of
``condensation'' in probability space, which is analogous to
the condensation phase in stochastic mass transport where a macroscopically
large mass forms at a single site on a lattice \cite{Majumdar09}.
In a stochastic process, the analogous phenomenon is the emergence
of one or more jumps responsible for a macroscopic fraction of the 
total displacement after many steps.
Condensation occurs in heavy-tailed processes, but can also 
occur even in light-tailed processes in the presence of multiple constraints
(e.g., conditioning on the values of both the total sum and the sum of squared
steps) \cite{Szavits14, Szavits14b}.
In another example, a double transition to the condensed state occurs
in the run-and-tumble particle \cite{Gradenigo19}.
Our results augment this prior work by providing a length scale defining
the crossover to the large-jump regime, which is based on the analysis of
conditioned probabilities.

\section{Definitions}

The continuous-time $\alpha$-stable Lévy processes are specified in terms of 
the characteristic function 
$\langle e^{ikx(t)}\rangle=e^{-t\sigma^\alpha |k|^\alpha}$ at time $t$, 
provided $x(0)=0$ \cite{jacobs10}; 
the Fourier transform yields the probability density 
$f_\alpha(x;t)$ for $x(t)$, 
thus being ``stable'' under iterated convolutions.
For simplicity we will only consider symmetric stable processes. 
Also, $\sigma$ is a width-scaling parameter, and $\alpha\in(0,2]$
characterizes the long tails of the densities.
The case $\alpha=2$ is Gaussian, while $\alpha<2$ densities
have heavy, power-law tails scaling as $|x|^{-(1+\alpha)}$.
The variance diverges for $\alpha<2$ and 
the mean absolute deviation diverges
for $\alpha\leq 1$.
The power-law tails are responsible for jump
discontinuities in the stochastic
evolution that are absent in the Gaussian case.
To be precise about terminology, 
we will refer to these jump discontinuities as ``jumps,''
while instead using ``steps'' or ``displacements'' to refer
to the change in state over a finite time interval.

\begin{figure}[t]
  \center
\includegraphics[width=\linewidth]{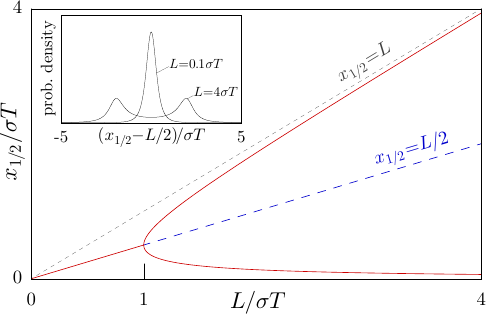}  \caption{
    Bifurcation diagram showing maxima (solid, red) and minima (dashed, blue) of the conditioned density~(\ref{levybridgemiddensity}) for $\alpha=1$.
    Inset: conditioned density before and after the bifurcation.
  \label{fig:cauchy_bifurcation}
}
\end{figure}

\begin{figure*}[t]
\center
\includegraphics[width=\textwidth]{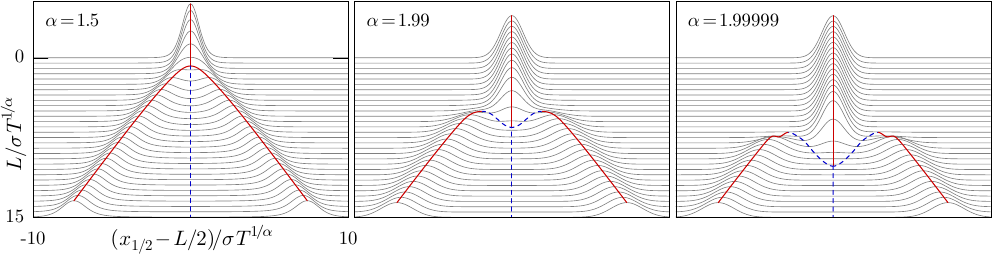}\caption{
  Variation of the midstep density (\ref{levybridgemiddensity}) with Lévy index
  $\alpha$ and arrival point $L$. Curves highlighting
  maxima (solid, red)
  and minima (dashed, blue) are superimposed. 
  \label{fig:waterfall3}
}
\end{figure*}

\section{Bifurcation Length}

\subsection{Lévy bridges}

In a Lévy bridge, the arrival point is specified as $x(T)=L$ for 
some arrival time $T>0$. Then the intermediate position $\xmid\!:=\!x(T/2)$
has the conditional density (``midpoint density'')
\begin{equation}
  {f}_\alpha(\xmid;T/2|x\!=\!L;T)=
  \frac{{f}_\alpha(\xmid;T/2)\,{f}_\alpha(L\!-\!\xmid;T/2)}{{f}_\alpha(L;T)}
  \label{levybridgemiddensity}
\end{equation}
in terms of the unconditioned density $f_\alpha(x;t)$.
Once $x(T/2)$ is sampled, the bridge is effectively bisected
into two bridges, and the midpoint-sampling process may be iterated
to sample the Lévy bridge to any desired time resolution.
For $\alpha=2$ the midpoint density retains the same Gaussian form
as the unconditioned density, but the conditioned and unconditioned
forms differ for any $\alpha<2$.

The Cauchy ($\alpha=1$) case is a good
example of what happens for $\alpha<2$.
For $L<\sigma T$,
this distribution has a single peak at $\xmid = L/2$, which has
a seemingly intuitive interpretation: if a particle 
travels from $x=0$ to $L$ in time $T$, the most probable 
intermediate position at $T/2$ is $L/2$. However, this
intuition breaks down at the special
arrival point $L_\mathrm{b}=\sigma T$, beyond which the
midpoint density becomes bimodal, and the single maximum
bifurcates into a pair at $\xmid=[L \pm (L^{2} - \sigma^2 T^2)^{1/2}]/2$
(Fig.~\ref{fig:cauchy_bifurcation}). 
For $L\gg L_\mathrm{b}$ the peaks are well separated, with maxima
approaching asymptotes $\xmid\sim 0,L$. In this case,
the interpretation of the midpoint changes: 
the large final displacement $L$ tends to break down into 
one large step of order $L$ and one small step,
rather than two steps roughly equal to $L/2$.
A bridge with sufficiently large overall 
transition length $L$ will tend to 
maintain this as a single jump discontinuity.

Similar structural changes in the midpoint density occur 
for all $\alpha<2$.  
Figure~\ref{fig:waterfall3} shows typical possibilities
of how the bifurcation occurs as $L$ increases. 
For $\alpha=1.5$ there is a pitchfork bifurcation
\cite{*[{The term ``bifurcation'' here refers to the behavior of the
maxima of the probability densities, behavior analogous to the bifurcation
of the stable points in a quartic potential.  This bifurcation terminology
has been applied to a probability density in the same sense by }] [{ in an economic stochastic-process model; by contrast, the bifurcation we discuss arises naturally from the conditioned stable Lévy densities themselves.}] chiarella91},
as in the Cauchy case, where two maxima and a minimum
are created from a single maximum.
However, closer to the Gaussian limit
($\alpha=1.99$ and $\alpha=1.99999$), the structure
is more complicated: first, a pair of
side peaks is born via tangent bifurcations; second,
the side peaks grow to match the central peak in height; and third,
a central minimum forms in a reverse-pitchfork bifurcation.
For any $\alpha$ the end results are the same:  a unimodal density transforms into
a bimodal density with well separated peaks.

\subsection{Variation with \boldmath$\alpha$} 
\begin{figure}[t]
  \begin{center}
\includegraphics[width=\linewidth]{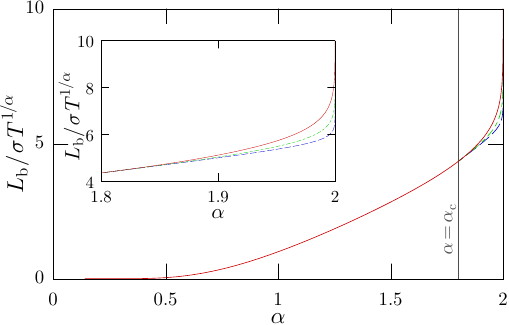}  \end{center}
  \vspace{-5mm}
  \caption{
    Variation of boundaries between ``small'' and ``large'' 
    steps with $\alpha$.  Curves indicate  
    bifurcation lengths $L_\mathrm{b}$ for which the center of the midstep 
    density has vanishing curvature (red, solid),
    the half-step distribution develops side peaks (blue, dashed), and
    side peaks are equal in height to the center peak (green, dot-dashed).
    Inset: magnified view for $\alpha>\alpha_\mathrm{c}$. 
    \label{fig:crit_L_vs_a}
  }
\end{figure}

An obvious characterization of the bifurcation length $L_\mathrm{b}$
is the value of $L$ for which the curvature of the midpoint density 
(\ref{levybridgemiddensity}) at $\xmid=L/2$ changes sign
(Fig.~\ref{fig:crit_L_vs_a}).
However, for $\alpha$ above a critical value $\alpha_\mathrm{c}$, as
we have seen, the midpoint density does not exhibit a simple
bifurcation to a bimodal density; 
rather, there are three distinct transitions. 
[The  critical value $\alpha_\mathrm{c}\approx 1.7999233$
occurs when the fourth derivative 
of the midpoint density (\ref{levybridgemiddensity}) vanishes at $\xmid=L/2$
(in addition to the vanishing of the second derivative, which already 
defines $L_\mathrm{b}$).]
All three bifurcation lengths are shown in 
Fig.~\ref{fig:crit_L_vs_a} for $\alpha>\alpha_\mathrm{c}$.
They all
usefully characterize the structural changes of
the distribution,
though in practice the particular choice of $L_\mathrm{b}$ is
not too important---as we will see, the 
transition between ``short'' and ``long'' displacements is not sharp.
(We use the curvature-change criterion except where noted.)

Figure~\ref{fig:crit_L_vs_a} also shows the transition
away from power-law tails in the limit $\alpha\longrightarrow 2$.
The bifurcation length diverges in this limit, so that
for the Gaussian $(\alpha=2)$ case, any final step $L$ is
a ``short step.''
The nature of this divergence may be analyzed
using the asymptotic density
$f_\alpha(x;t=1)\sim f_2(x;1) + \delta |x|^{\delta-3}$, 
valid for large $|x|$ and small $\delta:=2-\alpha$
\cite{nagaev89}.
One can show that $L_\mathrm{b}$
(defined by the curvature-sign-change criterion)
diverges as
$L_\mathrm{b}\sim [-4\sigma^2 T\log(\pi\delta^2/2)]^{1/2}$.
Numerically, $L_\mathrm{b}$ seems to diverge similarly according
to the other criteria as well.
Thus, even very close to the Gaussian limit $\alpha=2$,
$L_\mathrm{b}$ remains relatively small (cf.\  Fig.~\ref{fig:waterfall3}, third panel).

\subsection{Conditioned sampling}
\begin{figure}[tb]
  \center
\includegraphics[width=\linewidth]{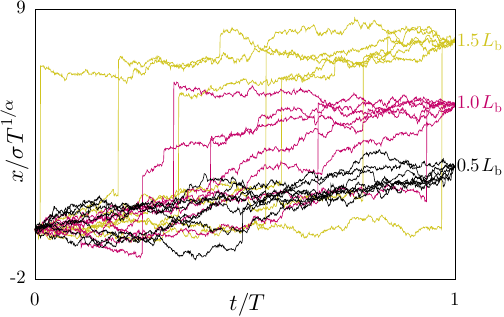}  \caption{
    Typical sample paths of Lévy bridges for $\alpha=1.9$,
    illustrating the qualitative transition with $L$.
    Each path was generated through 10 recursive subsamplings 
    from the midstep distribution (\ref{levybridgemiddensity}).
    \label{fig:bridges}
}
\end{figure}

As noted above, when sampling the intermediate state of a L\'evy bridge for $L>L_\mathrm{b}$, a jump of order $L$ likely persists. 
Upon further recursive subsampling of the bridge's intermediate
states, this behavior locks in: $L_\mathrm{b}$ is
effectively smaller when sampling sub-bridges on progressively smaller time intervals, so that
the substep length $L$ tends to exceed  $L_\mathrm{b}$ by an
ever increasing margin, making it progressively less likely
to be split into smaller jumps. 
Figure~\ref{fig:bridges} illustrates this:
for $L=1.5L_\mathrm{b}$ there is typically a single long
step that persists to high temporal resolution. By contrast,
for $L=0.5 L_\mathrm{b}$, the overall displacement has decomposed
into many small steps, with an appearance
resembling Brownian motion. 
The intermediate case $L=L_\mathrm{b}$ exhibits both behaviors.

This behavior under conditioned subsampling shows that the bifurcation length $L_\mathrm{b}$ yields an innate notion of
large steps of an $\alpha$-stable process. Specifically, an observed final displacement $|x(T)|\gg L_\mathrm{b}$ most likely corresponds to a single, similarly large jump discontinuity, even if the detailed evolution up to the final time $T$ is not known.  Meanwhile, a smaller final displacement  $|x(T)|\lesssim L_\mathrm{b}$ is more likely to be a composite event comprising multiple smaller jumps. 
This latter conclusion can be understood from the tails of the conditioned density (\ref{levybridgemiddensity}),
which scale as $|x|^{-2(1+\alpha)}$, which are relatively short compared to the
$|x|^{-(1+\alpha)}$ tails of the step density $f_\alpha(x;t)$.
This is a powerful qualitative inference based only on the endpoints of the process;
it is useful in 
problems of interpolation of a stochastic process between observations 
(e.g., animal movement \cite{horne07} and kriging \cite{stein99}),
if the underlying process is heavy-tailed.
 Additionally, this provides a means for inferring whether
 a rare, significant event occurred between observations. 
 Such criteria are important for the  analysis of statistical extremes \cite{beirlant2004}
 and for specific problems like detecting market crashes \cite{schluter08}. 

A salient feature of stable Lévy processes is scale-invariance. 
So how is it possible to have an intrinsic scale $L_\mathrm{b}$?
Scale invariance is best seen in the Lévy--Khintchine representation \cite{cont04, gardiner09, jacobs10}, 
where symmetric stable processes have pure power-law jump-rate densities
$\sin(\pi \alpha/2)\Gamma(1+\alpha)\,\sigma^\alpha/(\pi|\Delta x|^{1+\alpha})$.
In some sense, then, any scale based solely on the step distribution (width at half maximum, etc.)
is inherently nonsensical. 
However, conditioning introduces a timescale $T$, which 
induces a length scale---one that can only be understood through the
variable structure of the conditioned density (\ref{levybridgemiddensity}).
Importantly, this scale differs from the well known length scale $\sigma T^{1/\alpha}$ \cite{Vezzani19, Burioni20}.
This distinction defines an intuitive notion of ``long'' displacements 
that captures
how, visually and intuitively, the large-scale structure of stable Lévy
processes seem similar to Gaussian processes punctuated by discrete
jump discontinuities (Fig.~\ref{fig:jumps_highlight}). 
Mathematically, this similarity is not obvious:
Stable Lévy processes with $\alpha<2$ have a dense set
of discontinuities, whereas Gaussian process are 
continuous (almost surely).

\section{Applications}

\subsection{Stretched Lévy bridges}
In the Gauss\-ian case, one important representation of the Brownian
bridge is \cite{karatzas1991} 
\begin{equation}
  W(t) = B(t) + \frac{t}{T}\big[W(T)-B(T)\big],
  \label{WB}
\end{equation}
where $W(t)$ is a Wiener process (unconditioned 
Lévy process with $\alpha=2$, $\smash{\sigma=1/\sqrt{2}}$),
and $B(t)$ is a Brownian bridge 
[Wiener process conditioned to have a fixed arrival $B(T)$].
Intuitively, in the ``standard bridge'' case $B(T)=0$, 
the second term is the ballistic trajectory from $0$ to $W(T)$,
while $B(t)$ comprises the random fluctuations.
This representation, when interpreted as an expression for
$B(t)$ in terms of $W(t)$ and the ballistic motion,
provides a simple way to
simulate Brownian bridges using any
Wiener-process algorithm.  
Naively, it seems like this representation should be valid for $\alpha<2$
stable processes:
Dividing the evolution into time steps $\Delta t$, the increments
of the stable process and bridge are of order $\Delta t^{1/\alpha}$,
while the ballistic correction is of order $\Delta t$.
The ballistic component is thus of order  $\Delta t^{1-1/\alpha}$ relative
to the Lévy-process steps, and thus should be negligible as $\Delta t\longrightarrow 0$
provided $\alpha>1$.
In the Gaussian case this heuristic
argument is correct, and the representation (\ref{WB}) is
valid---any ballistic ``stretch'' does not affect the Gaussian
statistics in the continuum limit. It fails, however, for $\alpha<2$:
if $x(T)$ corresponds to a sufficiently large final displacement, 
then the stretch is excessive, and the resulting ``bridges''
produce erroneous results in simulations.
(Reference~\cite{knight96} noted this 
inequivalence between stretched and conditioned bridges 
\cite{*[{Another inequivalent
representation was studied by }] [{}] janicki94}.) 

\begin{figure}[tb]
  \center
\includegraphics[width=\linewidth]{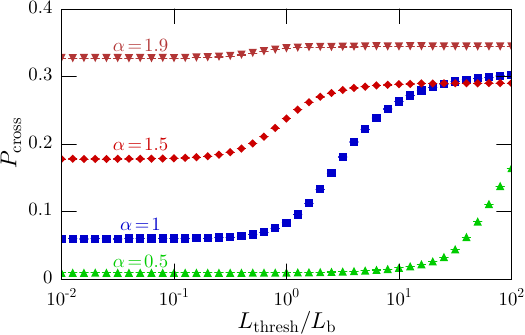}  \caption{
    Simulated probability for Lévy bridges generated
    via Eq.~(\ref{WB}) to cross a boundary at
    $d=\sigma T^{1/\alpha}$, as the rejection
    threshold $L_\mathrm{thresh}$ varies.  
  \label{fig:first_cross_stretch}
}
\end{figure}

Since we have a large-step criterion, it
is possible to deal with excessive stretches.
The fix is to define a threshold 
$L_\mathrm{thresh}$, and an unconditioned Lévy sample path is only
stretched as in Eq.~(\ref{WB}) if its final point $L=x(T)$
is within $L_\mathrm{thresh}$ of the bridge's arrival point.
Otherwise, it is rejected and other paths attempted until a 
bridge is successfully generated.
The $\alpha$-dependent bifurcation length $L_\mathrm{b}$ from
Fig.~\ref{fig:crit_L_vs_a}
marks a scale $L_\mathrm{thresh}$ below which the stretching
algorithm should yield an accurate set of Lévy bridges.
A test of this algorithm, computing the probability $P_\mathrm{cross}$
for Lévy bridges (with $L=0$) 
to cross a boundary at $d=\sigma T^{1/\alpha}$ before
time $T$, illustrates this transition
(Fig.~\ref{fig:first_cross_stretch}) \footnote{Simulations used $\Delta t=10^{-5}T$,  averaging over $10^7$ paths.}. 
In particular, the
simulated $P_\mathrm{cross}$ rapidly becomes accurate when $L_\mathrm{thresh}$
decreases below $L_\mathrm{b}$ 
(the bridge construction is exact in the limit $L_\mathrm{thresh}\longrightarrow 0$).
As a practical bridge-generation method, this is much more
efficient than using $\sigma\Delta t^{1/\alpha}$
(the smallest natural length scale) for $L_\mathrm{thresh}$.

A particularly interesting feature in Fig.~\ref{fig:first_cross_stretch}
is that $P_\mathrm{cross}=0.9\%$ is so small for the case
$\alpha=0.5$.
(By contrast, $P_\mathrm{cross}=31.5\%$ in the unconditioned case.)
The surprise here is that the smallest-$\alpha$ case has the strongest tendency towards
large jumps---intuitively, the best ``mobility''---and yet has the smallest boundary-crossing probability.
However, conditioning on $L=0$ also
conditions away the tendency to have extreme jumps (and thus to easily cross the boundary), precisely because an extreme
jump is suppressed by the requirement of a compensating (and correspondingly rare) jump to return to the final target state.

\begin{figure}[tb]
  \center
\includegraphics[width=\linewidth]{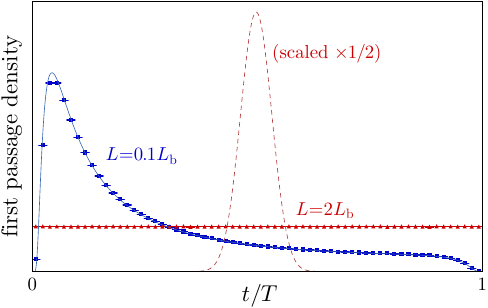}  \caption{
    Simulated conditioned first passage time distributions 
    for $\alpha=1.99999$ and
    $d=L/2$ are shown for $L=0.1 L_\mathrm{b}$ 
    (blue/squares) and  $L=2L_\mathrm{b}$ (red/triangles). 
    Exact densities for $\alpha=2$ \cite{borodin02}
    for the same $L$ values are shown for comparison in each
    case (blue, solid and red, dashed, respectively).
  \label{fig:fpts}
}
\end{figure}

\subsection{Conditioned first passage}\label{section:conditionedFP}
First-passage times, defined here as the first time a process $x(t)$ exceeds a boundary $d$, are of broad importance \cite{redner01}.
They are especially interesting for Lévy
processes due to the universal Sparre Andersen scaling \cite{zumofen95,chechkin03,klafter11}, where the tail of the 
first-passage-time distribution is $\alpha$-independent. 
However, as we have seen, conditioned 
Lévy bridges have a particularly sensitive transition as 
$\alpha\longrightarrow 2$, a pattern that continues
for first-passage times.

An intuitive picture of the conditioned first-passage time follows from the qualitative appearance of the sample paths for $L=1.5L_{\mathrm{b}}$ in Fig.~\ref{fig:bridges}.  A dominant jump is consistently present among the paths, but not at any particular time. This can be regarded as an outcome of recursively sampling the midpoint density (\ref{levybridgemiddensity}). For $L\gg L_\mathrm{b}$, a large step likely persists under sampling iterations, but due to the symmetry of the midstep distribution, the large step is equally likely to be associated with any time subinterval. Since the first-passage time is likely due to the dominant jump, the first-passage time should be uniformly distributed. 
Figure~\ref{fig:fpts} confirms this intuition with simulations of the first passage density \footnote{Simulations averaged $10^7$ paths, with $\Delta t=2^{-14}T$.}. For $L=2L_\mathrm{b}$ the first passage density is indeed uniform.  A small change from $\alpha=1.99999$ to the Gaussian case yields a remarkably different distribution: approximately Gaussian, centered at $t\approx T/2$.  The Gaussian result follows intuitively from 
Eq.~(\ref{WB}), since the most likely bridges in this regime are concentrated around the ballistic path to the endpoint.

For a smaller overall displacement ($L=0.1L_\mathrm{b}$), 
the first-passage-time densities in the
$\alpha=1.99999$ and Gaussian cases match closely. 
This is consistent with the
observation that for $L\ll L_\mathrm{b}$, the conditioned Lévy bridges are qualitatively similar to Brownian bridges.
Nevertheless, the rare but important extreme jumps
generate remarkably non-Gaussian behavior, even close
to the Gaussian limit.

\subsection{Unconditioned first escape}

\begin{figure}[tb]
  \center
\includegraphics[width=\linewidth]{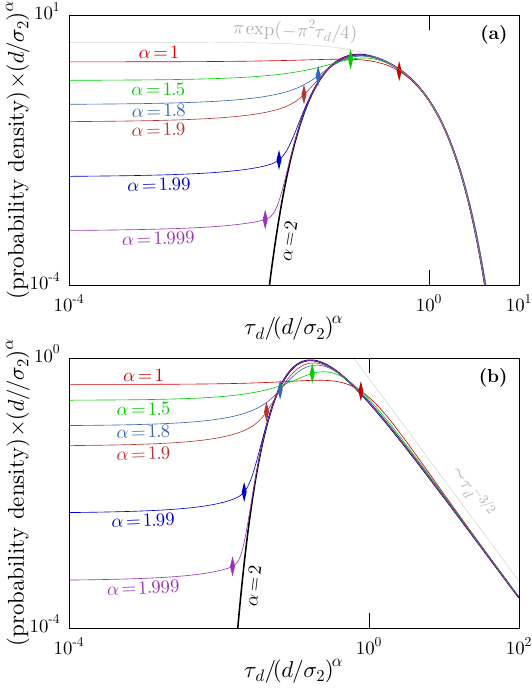}  \caption{
    (a)
    Simulated first-escape-time distributions from $(-d,d)$
    for verious Lévy processes starting at $x=0$, showing 
    long-time behavior that matches Gaussian ($\alpha=2$) behavior,
    but diverging short-time behavior.
    (b) Simulated first-passage distributions [i.e., first escape from
    $(-\infty,d)$] show similar behavior.
    The bifurcation time $T_\mathrm{b}$ is marked by a diamond
    symbol in each case.
    The width-scaling parameter $\sigma$ is chosen for each $\alpha$
    to obtain matching asymptotic behavior, as discussed in the main 
    text.
  \label{fig:fetsfpts}
}
\end{figure}

The major theme of this paper has been the conditioned evolution
of stochastic processes. However, the reasoning we have used thus far
is useful in studying unconditioned evolution as well.
As a common example, consider the first-escape time of a stable Lévy process
starting at $x=0$
from the interval $(-d,d)$, $\tau_d:=\inf\{t:|x(t)|\geq d\}$.
The simulated probability densities for the escape
time for various values of $\alpha$ are shown in 
Fig.~\ref{fig:fetsfpts}(a).
The striking feature of this set of distributions is the universal,
$\alpha$-independent asymptotic behavior at long times.
Each of the probability distributions splits away from the 
Gaussian ($\alpha=2$) distribution at a point that is different for
each value of $\alpha$.

To understand the behavior here, first note that the
portion of the escape-time distribution to the left of any particular
time $\tau$ acts as a conditioned density, because it refers
only to the subset of trajectories that has escaped by time $\tau$.
This implicit conditioned behavior allows us to apply our results for
conditioned processes to this simple escape-time problem.
The bifurcation length $L_\mathrm{b}$ is of particular utility
here.
Recall that a long jump must have occurred in an escape by time
$\tau$ if $d>L_\mathrm{b}=\tilde{L}_\mathrm{b}\sigma \tau^{1/\alpha}$,
where
$\tilde{L}_\mathrm{b}$ is the value of the bifurcation length
$L_\mathrm{b}$ given by setting $\sigma=T=1$.
Rearranging this expression, we can define the bifurcation timescale
$T_\mathrm{b}$ such that an escape by time $\tau$ must have involved
a long jump if
$\tau<T_\mathrm{b}:=(d/\sigma)^\alpha/\tilde{L}_\mathrm{b}^{\,\alpha}$.
This bifurcation time is marked as a diamond on each distribution
in Fig.~\ref{fig:fetsfpts}(a), and it evidently
marks the timescale where the escape-time distribution
for each stable Lévy case splits away from the Gaussian limit.
For large escape times $\tau_\mathrm{d}\gg T_\mathrm{b}$,
an extreme jump is unlikely, and any Lévy process behaves basically as a 
Gaussian random walk.
On the other hand, for small escape times
$\tau_\mathrm{d}\ll T_\mathrm{b}$,
a single large jump is the most likely scenario.
In this case the reasoning of Sec.~\ref{section:conditionedFP}
applies, and the dominant jump is equally likely to occur
at any time below a fixed $\tau\ll T_\mathrm{b}$.
In the escape-time distributions, this behavior appears as an asymptotically
constant behavior of the distribution as $\tau_d\longrightarrow 0$
(where in the Gaussian case, the probability density vanishes here).
This constant value of the density at small escape times decreases with
increasing $\alpha$, as expected because the probability of an extreme
jump also decreases (owing to the less-fat tails).
The universality of the long-time asymptotic tail here is thus another example
of the intuition we mentioned above that heavy-tailed processes 
resemble Gaussian processes between occurrences
of rare, extreme events.

Figure~\ref{fig:fetsfpts}(b) shows the analogous behavior for
(unconditioned) first-passage
densities $\tau_d:=\inf\{t:x(t)\geq d\}$, corresponding to escape from
the interval $(-\infty,d)$.
The division between Gaussian-like and extreme-event behavior is also apparent
here---the main difference is the form of the asymptotic tail, which has
the characteristic Sparre Andersen scaling of $\smash{\sim\!\tau_d^{\,-3/2}}$.

At this point some brief comments clarifying the simulated distributions
in Fig.~\ref{fig:fetsfpts} are in order.
The distributions were computed by
numerical integration of the fractional diffusion equation \cite{Watanabe62}.
It is most sensible to compare distributions with
the same long-time asymptotic behavior, accomplished by
an appropriate choice for the
$\alpha$-dependent  width-scale parameter $\sigma_\alpha$.
For the first-escape problem, the asymptotic tail
is of the form $\smash{e^{-\sigma_\alpha^{\,\alpha}\lambda_1^{(\alpha)} t/d^\alpha}}$ \cite{Dybiec17},
where $\smash{\lambda_1^{(\alpha)}}$ is  
the smallest eigenvalue of the fractional Laplace operator
$(-\nabla^2)^{\alpha/2}$ on the bounded domain $[-1,1]$
(e.g., $\smash{\lambda_1^{(2)}=\pi^2/4}$ in the Gaussian limit).
The asymptotics thus
match across $\alpha$ via the choice
$\smash{\sigma_\alpha^{\,\alpha}/\sigma_2^{\,2}
  =\pi^2/4d^{2-\alpha}\lambda_1^{(\alpha)}}$,
using the Gaussian scale parameter $\sigma_2$ as a reference.
The asymptotic tail in the first-passage problem 
has the form 
$(d/\sigma)^{\alpha/2}/\alpha\sqrt{\pi}\Gamma(\alpha/2)\tau^{3/2}$
\cite{koren07}.
These asymptotics match for the choice
$\smash{\sigma_\alpha^{\,\alpha}/\sigma_2^{\,2}=
  1/d^{2-\alpha}\Gamma^2(1+\alpha/2)}$.

The bifurcation time (and thus length) mark a boundary
between Gaussian and extreme-event behavior
in a conceptually simple setting of the escape problem,
which is directly amenable to experimental observation.
For example, we have already mentioned that laser-cooled atoms
are an important prototype system for studying Lévy-type dynamics
\cite{marksteiner1996, katori1997, bardou02, sagi2012, kessler2012, barkai2014, afek2017, aghion2017}
(including the big-jump principle \cite{Vezzani19}).
A setup appropriate for the study of escape times is that of a
single laser-cooled atom monitored by a 
fluorescence-detection system \cite{katori1997, Alt02}.
An aperture for the fluorescence photodetector defines the region from which the
atom is to escape; the time that it takes to observe a sharp drop
in the atomic fluorescence after the release of the atom (from its initially
prepared position) is a measure of the escape time.
A more detailed discussion of the Lévy behavior of laser-cooled atoms
as well as typical parameters for an experimental realization are included
in the Supplemental Material \cite{supplement}.
Of course, beyond laser-cooled atoms, this first-escape behavior should
be observable in essentially any stochastic physical system that is accurately
modeled by stable Lévy processes.

\section{Summary}

We have discussed the conditioned evolution of $\alpha$-stable Lévy 
processes as a prototype for extreme events.
The knowledge of a particular final state turns out to retrodict
whether an extreme event occurred along the way.  
This conclusion follows from an analysis of the
conditioned densities, which change form as the final displacement passes
a threshold, the bifurcation length.
The analysis here has applications to the construction of Lévy bridges,
the qualitative understanding of conditioned first-passage dynamics,
and the understanding of unconditioned first-escape problems.
We have also pointed out how the manifestation of the bifurcation length
in the first-escape problem can be studied experimentally
with laser-cooled atoms.

\begin{acknowledgments}
We gratefully acknowledge helpful discussions with Steven van Enk.
This work was supported by the NSF (PHY-1505118) and NVIDIA Corporation.
\end{acknowledgments}

\bibliographystyle{apsrev4-2}
\bibliography{bifurc}

\begin{thebibliography}{66}
\bibitem [{sup()}]{supplement}%
  \BibitemOpen
  \href@noop {} {\bibinfo {title} {See {S}upplemental {M}aterial for a
  discussion of the observability of the effects described in this paper in
  experiments with laser-cooled atoms, which includes {R}efs.~\cite{Dalibard89,
  Jonsella06, Svensson08, Lutz03}.}}\BibitemShut {Stop}%
\end{thebibliography}

\clearpage
\pagestyle{empty}
\widetext
\includepdf[pages=1]{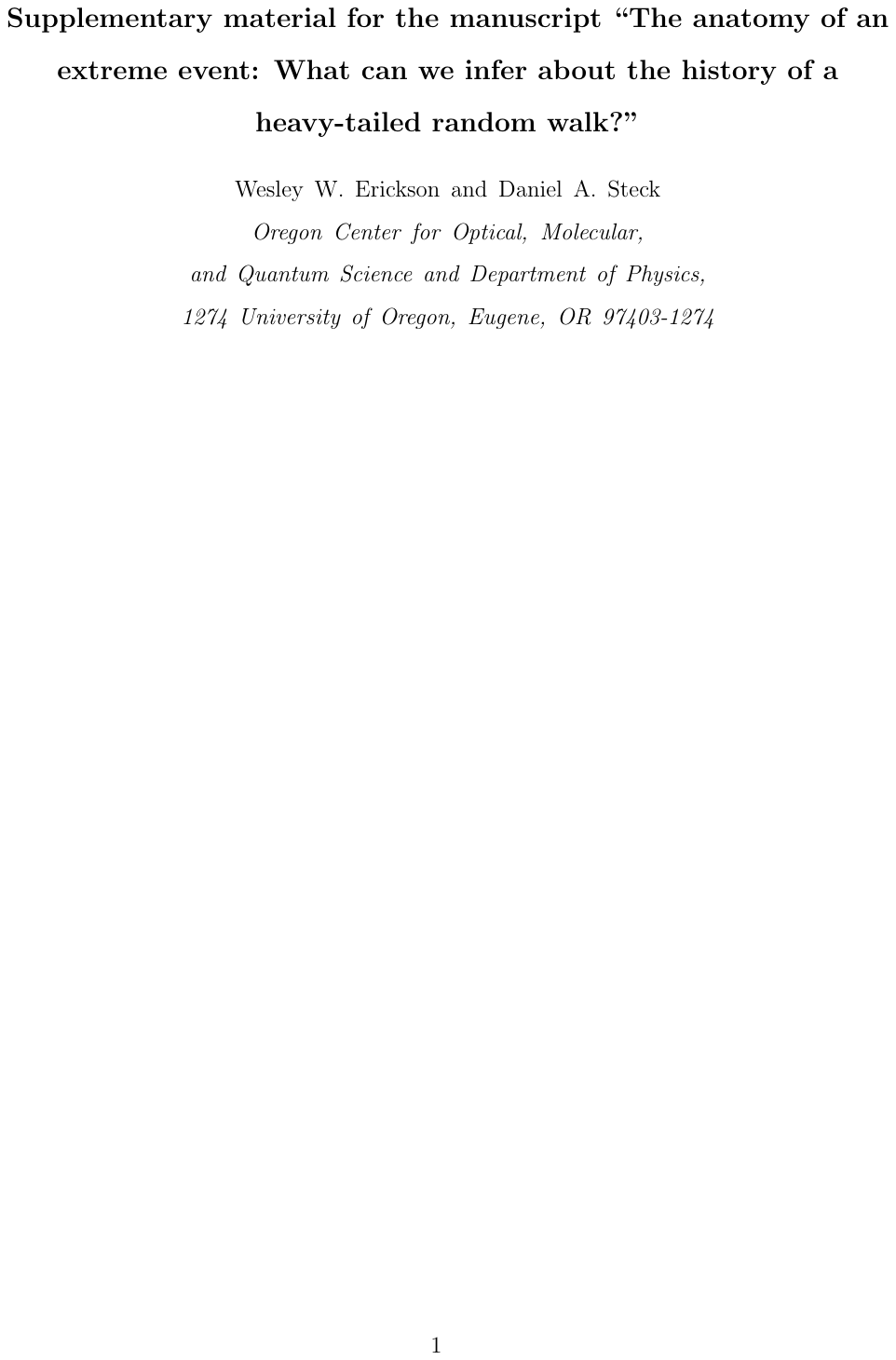}
\includepdf[pages=2]{supplement.pdf}
\includepdf[pages=3]{supplement.pdf}
\includepdf[pages=4]{supplement.pdf}
\includepdf[pages=5]{supplement.pdf}
\includepdf[pages=6]{supplement.pdf}
\includepdf[pages=7]{supplement.pdf}
\includepdf[pages=8]{supplement.pdf}

\end{document}